%% file: Newton_for_Networks.tex
\input TEXSHOP_macros_new.tex

\def\section#1{\goodbreak\vskip 3pc plus 6pt minus 3pt\leftskip=-2pc
   \global\advance\sectnum by 1\eqnumber=1\subsectnum=0%
\global\examplnumber=1\figrnumber=1\propnumber=1\defnumber=1\lemmanumber=1\assumptionnumber=1 \conditionnumber =1%
   \line{\hfuzz=1pc{\hbox to 3pc{\bf 
   \vtop{\hfuzz=1pc\hsize=38pc\hyphenpenalty=10000\noindent\uppercase{\the\sectnum.\quad #1}}\hss}}
			\hfill}
			\leftskip=0pc\nobreak\tenf
			\vskip 1pc plus 4pt minus 2pt\noindent\ignorespaces}
\def\subsection#1{\noindent\leftskip=0pc\tenf
   \goodbreak\vskip 1pc plus 4pt minus 2pt
               \global\advance\subsectnum by 1
   \line{\hfuzz=1pc{\hbox to 3pc{\bf \the\sectnum.\the\subsectnum.
   \vtop{\hfuzz=1pc\hsize=38pc\hyphenpenalty=10000\noindent{\bf #1}}\hss}}
                        \hfill}
   \leftskip=0pc\nobreak\tenf
                        \vskip 1pc plus 4pt minus 2pt\nobreak\noindent\ignorespaces}



\input miniltx

\ifx\pdfoutput\undefined
  \def\Gin@driver{dvips.def} 
\else
  \def\Gin@driver{pdftex.def} 
\fi

\input graphicx.sty
\resetatcatcode

\long\def\fig#1#2#3{\vbox{\vskip1pc\vskip#1
\prevdepth=12pt \baselineskip=12pt
\vskip1pc
\hbox to\hsize{\hfill\vtop{\hsize=30pc\noindent{\eightbf Figure #2\ }
{\eightpoint#3}}\hfill}}}

\def\show#1{}

\rightheadline{\botmark}

\pageno=1

\rightheadline{\botmark}

\pn {\bf April 2011}\hfill{\bf Report LIDS - 2866}
\bigskip \bigskip\bigskip

\bigskip

\def\longpapertitle#1#2#3{{\bf \centerline{\helbigb
{#1}}}\medskip{\bf \centerline{\helbigb
{#2}}}\bigskip{\bf \centerline{
{#3}}}\bigskip}

\vskip-3pc

\longpapertitle{Centralized and Distributed Newton Methods for}{
Network Optimization and Extensions}{ {Dimitri P.\ Bertsekas\footnote{\dag}{\ninepoint  Dimitri Bertsekas is with the Dept.\ of Electr.\ Engineering and
Comp.\ Science, M.I.T., Cambridge, Mass., 02139. His research was supported by the Air Force Grant FA9550-10-1-0412.} }}

\centerline{\bf Abstract}

We consider Newton methods for common types of single commodity and multi-commodity network flow problems. Despite the potentially very large dimension of the problem, they can be implemented using the conjugate gradient method  and low-dimensional network operations, as shown nearly thirty years ago. We revisit these methods, compare them to more recent proposals, and describe how they can be implemented in a distributed computing system. We also discuss generalizations, including the treatment of arc gains, linear side constraints, and related special structures.

\vskip1pc

\section{Introduction}

\vskip-0pc

\pn A common type of nonlinear network flow optimization involves a graph with a set of directed arcs ${\cal A}$, and a set of paths ${\cal P}$. Each path $p$ consists of a set of arcs $A_p\subset {\cal A}$. A set of path flows $x=\{x_p\mid p\in{\cal P}\}$ produces a flow at each arc $a\in {\cal A}$, 
$$f_a=\sum_{\{p\mid a\in A_p\}}x_p.$$
The objective is to find $x$ that minimizes
$$F(x)=\sum_{p\in{\cal P}}R_p(x_p)+\sum_{a\in{\cal A}}D_a(f_a),\xdef\costfn{\lab}\eqnum\show{oneo}$$
where $D_a$ and $R_p$ are twice continuously differentiable functions. Often either $R_p=0$ or else $R_p$ encodes a penalty for a range of values of $x_p$ (such as small nonnegative values), while $D_a$ encodes a penalty for a range of values of $f_a$ (such as values close to a ``capacity" of arc $a$). The minimization may be subject to some constraints on $x$, such as upper and lower bounds on each $x_p$, supply/demand constraints (sum of path flows with the same origin and destination must take a given value), side constraints etc. Sometimes capacity constraints on the arc flows are imposed explicitly rather than through the functions $D_a$. This is a standard framework, discussed for example in the author's works ([BeT89], Section 7.6, [BeG92], Section 5.7.3, [Ber98], pp.\ 391-398), and in many other works, including the survey [FlH95], which contains extensive references, and textbooks such as [Roc84] and [Pat98].

There are several methods for solving this problem. Among first-order gradient-based methods, we note the conditional gradient/Frank-Wolfe method [FGK73], [FlH95], gradient projection methods with and without diagonal scaling [Gal77], [Ber80], [BGG84], [BeG92], and simplicial decomposition/inner approximation methods [CaG74], [Hol74], [Hoh77], [HLV87], [Pat98], [BeY10].  Gradient projection methods are well-suited for distributed implementation, and several such methods have been proposed for routing and flow control,  dating to the early days of data networking  (see [Ber80], [GaG80], [BGG84], [TsB86]). Similar methods have also been proposed for traffic equilibrium analysis, where the same types of optimization problems arise (see the survey [FlH95] and the book [Pat98]). 

In this note, we focus on second order/Newton-type methods, which have also been proposed since the early days of data networking, but have received less attention, possibly because of the emphasis on simple methods in large-scale optimal network flow problems. A common concern with Newton-type methods is that they require the calculation of the gradient and Hessian of $F$, and the solution of some large-scale quadratic optimization problem, either unconstrained or constrained. However, what has been sometimes overlooked is that one may exploit special problem structure to calculate the Newton direction without forming or inverting the large-dimensional Hessian matrix of $F$. This idea is quite common in the numerical analysis and interior point method literature, where the Newton direction is often calculated using essentially iterative methods, such as conjugate direction, splitting, and Krylov subspace methods. We will discuss primarily the unconstrained case and refer to the literature for various ways to handle constraints, either directly or through the use of penalty and augmented Lagrangian functions. Diagonal second-order preconditioning may also be used to make the method more effective, and to facilitate the choice of stepsize.

Our point of departure is a fact known since the papers by Bertsekas and Gafni [BeG83], [GaB84], which proposed Newton-type methods for various types of constrained optimization problems, including multicommodity flow problems that involve nonnegativity constraints on $x$ and supply/demand constraints. Reference [BeG83] showed that for any $x$, the pure Newton direction $y_N$, which satisfies
$$\gr^2F(x)\,y_N=-\gr F(x),\xdef\newton{\lab}\eqnum\show{oneo}$$
can be conveniently calculated by using the conjugate gradient method, and graph operations that do not require the explicit formation and storage of $\gr^2F(x)$ or its inverse. Reference [GaB84] embedded this idea within a broader class of two-metric gradient projection algorithms, and provided computational results. In practice, the Newton direction should be approximated by only a limited number of conjugate gradient iterations, thereby obtaining an approximate Newton direction, which however is still a descent direction, regardless of how many conjugate gradient iterations are performed. Diagonal second-order preconditioning may be used to make the method more effective, and to facilitate the choice of stepsize.

The computational complexity per iteration of the conjugate gradient method is $O(T)$, where $T$ is the total number of arcs traversed by the paths in ${\cal P}$. The complexity of computing the Newton direction is thus $O(PT)$, where $P$ is the number of paths in ${\cal P}$ (assuming the worst case number $P$ for the number of conjugate gradient iterations to compute the Newton direction). If the ``average" number of arcs $T/P$ traversed by a path is much smaller than $P$, the conjugate gradient-based $O(PT)$ computation is much faster than the $O(P^3)$ computation required to obtain the Newton direction by direct matrix inversion.  If only  $m$ conjugate gradient iterations are performed to compute approximately the Newton direction, the total required number of computational operations  is $O(mT)$. These estimates should be compared with the computation required to compute the gradient of the cost function and to perform a diagonally scaled gradient iteration, which is $O(T)$.

In this note, we will also describe how the graph operations for conjugate-gradient based Newton direction calculation can be done in a distributed computing environment, where there is a processor assigned to each path and a processor assigned to each arc. Then the Newton direction can be calculated simply, with each processor updating, maintaining, and communicating to other processors just a few numbers.  Moreover the processor assigned to path $p$ computes the corresponding component $y_p$ of the Newton direction, and may execute (in synchrony with the other processors) its local portion of a pure Newton iteration. The processors may also collaborate to compute a stepsize guaranteeing descent at every iteration. Of course, it is possible to assign multiple paths and arcs to a single processor, with an attendant increase of the computation and communication load of the processor. Information exchange between the processors may be conveniently and accurately performed by using a spanning tree with a designated node serving as a synchronizer and leader for the distributed computation (see [BeT89], [BeG92]).

After describing the centralized and distributed methods for calculating the Newton step, we explore various generalizations, involving arc gains, side constraints, etc. Part of this discussion has not explicitly appeared in the literature so far, and we record it here for the purpose of making a connection with recent works on Newton-like methods for distributed routing and flow control in data networks [AtL00], [JOZ09], [WOJ10]. At the same time, we should note that most of what we describe, including the central idea, was essentially known in 1983. 

\vskip-1pc
\section{Centralized Computation of the Newton Step}

\pn We focus on a fixed path flow vector $x$ and to simplify notation, we suppress the arguments $x$, $x_p$, and $f_a$ from various function expressions. We show how to calculate the Newton step of Eq.\ \newton\ using the conjugate gradient method and graph-based operations. We write
$$g=\gr F(x),\qquad H=\gr^2 F(x),$$
and for each  path $p$, we denote by $g_p$ and $H_{pp}$ the corresponding component of $g$ and diagonal component of $H$.

We have
$$g_p=R_p'+\sum_{a\in A_p}D_a',\qquad p\in{\cal P},\xdef\gradcomp{\lab}\eqnum\show{oneo}$$
where $R_p'$ and $D_a'$ denote the first derivatives of $R_p$ and $D_a$,  calculated at the current arguments $x_p$ and $f_a$, respectively. 
Moreover 
$$H=\gr^2 R+E'\gr^2D E,\xdef\hescomp{\lab}\eqnum\show{oneo}$$
where $\gr^2R$ is the diagonal matrix with the second derivatives $R_p^{''}$ along the diagonal, $\gr^2D$ is the diagonal matrix with the second derivatives $D_a^{''}$ along the diagonal, and $E$ is the matrix with components
$$E_{ap}=\cases{1&if $a\in A_p$,\cr
0&if $a\notin A_p$.\cr}$$
Of particular interest are the diagonal elements of $H$, which are
$$H_{pp}=R_p^{''}+\sum_{a\in A_p}D_a^{''},\qquad p\in{\cal P}.\xdef\diaghescomp{\lab}\eqnum\show{oneo}$$
Note that the computation of $g_p$ and $H_{pp}$, collectively for all $p$, requires a total of $O(T)$ computational operations.

By exploiting the special structure of $g$ and $H$, we can perform some computations of interest using graph operations. These are:

\nitem{(a)} For each path $p$, we can calculate $g_p$ by accumulating first derivatives along $p$ [cf.\ Eq.\ \gradcomp]. Moreover, we can similarly calculate the diagonal elements $H_{pp}$  by accumulating second derivatives along $p$ [cf.\ Eq.\ \diaghescomp]. The number of required computational operations  is $O(T)$.

\nitem{(b)} For any vector $v=\{v_p\mid p\in{\cal P}\}$,  we can calculate the matrix-vector product
$$w=Hv$$
[cf.\ Eq.\ \hescomp] as the gradient of the function
$$\half v'Hv=\half \sum_{p\in {\cal P}}R_p^{''}v_p^2+ \half \sum_{a\in {\cal A}}D_a^{''}f^2_{a,v},$$
where 
$$f_{a,v}=\sum_{\{p\mid a\in A_p\}}v_p,\qquad a\in{\cal A}.\xdef\flowacum{\lab}\eqnum\show{oneo}$$
Thus, once we calculate $f_{a,v}$ for all ${a\in {\cal A}}$ by accumulating $v_p$ along path $p$ [an $O(T)$ computation], we  can obtain the $p$th component $w_p$ of $w=Hv$ as
$$w_p=R_p^{''}v_p+\sum_{a\in A_p}D_a^{''}f_{a,v},\qquad p\in{\cal P},\xdef\wpcomp{\lab}\eqnum\show{oneo}$$
using a similar computation to the gradient \gradcomp\ and the diagonal Hessian components \diaghescomp. Again the total number of required computational operations  is $O(T)$.

We now turn to the calculation of the Newton direction by viewing it as the solution of the quadratic problem
$$\min_y\,C(y)\;{\buildrel\rm def\over=}\;g'y+{1\over 2} y'Hy.\xdef\cycomp{\lab}\eqnum\show{oneo}$$
The idea is to solve this problem with the conjugate gradient method, which requires just inner products and matrix-vector products of the form $Hv$ that can be done using graph operations, without explicit formation of $H$, as we have seen. We assume that $H$ is positive definite so that the Newton direction is defined and the subsequent conjugate gradient algorithm has guaranteed termination, but this assumption is not essential for constrained versions of Newton-type methods, and for problems of practical interest, related conjugate-gradient based algorithms can be obtained (see e.g., [BeG83], [GaB84]).\footnote{\dag}{\ninepoint  If $H$ is positive semidefinite but singular, either the problem \cycomp\ has an infinite number of optimal solutions, in which case the conjugate gradient method will find one in a finite number of iterations (at most $P$ minus the dimension of the optimal solution set), or else problem \cycomp\ has no solution, in which case for some $k$, the stepsize $\a_k$ in the subsequent conjugate gradient method will become ``infinite"along a conjugate direction $p_k$ where $p_k'Hp_k=0$. Then by moving along the line $y_k+\a p_k$, we will be continuously decreasing the cost $C(y)$, and any vector of the form $y_k+\a p_k$ with $\a>0$ will be a direction of descent of $F(x)$ at $x$, which may be used in place of the Newton direction. Conjugate gradient methods are often implemented with such safeguards against singularity of $H$,  particularly in constrained contexts.} 

In particular, we start the iterative process with $y_0=0$ as the initial iterate for the minimum of the quadratic cost $C$, with the gradient of $C$ at $y=0$, which is $r_0=g$, and with $p_0=-g$, which we view as the first conjugate direction. Given the current iterate-gradient-conjugate direction triplet $(y_k,r_k,p_k)$, we generate the next iterate-gradient-conjugate direction triplet $(y_{k+1},r_{k+1},p_{k+1})$ by a conjugate gradient iteration, which has the following form (see e.g., [Lue84], [Ber99], Section 1.6):
$$\eqalign{y_{k+1}&=y_k+\a_kp_k,\qquad\hbox{where }\a_k={r_k'r_k\over p_k'Hp_k},\cr
r_{k+1}&=g+Hy_{k+1},\cr
p_{k+1}&=-r_{k+1}+\b_kp_k,\qquad\hbox{where }\b_k={r_{k+1}'r_{k+1}\over r_k'r_k}.\cr}$$
Here $\a_k$ is the stepsize that minimizes $C(y)$ over the line $\{y_k+\a p_k\mid a\in\re\}$; it can also be written equivalently (but less conveniently) as $\a_k=-{p_k'r_k/ p_k'Hp_k}$ (since $p_k=-r_k+\b_k p_{k-1}$ and $r_k$ is orthogonal to $p_{k-1},\ldots,p_0$, a basic property of the conjugate gradient method). The  matrix-vector products to be computed at each iteration are $Hp_k$ and $Hy_{k+1}$, and they can be obtained by using Eqs.\ \flowacum-\wpcomp, with $v=p_k$ and $v=y_{k+1}$, respectively.

According to well-known theory [Lue84], [Ber99], either $y_{k+1}$ minimizes $C(y)$ (i.e., $r_{k+1}=0$) and hence $y_{k+1}$ is equal to the Newton direction. Moreover (unless $g=0$) we have
$$C(y_k)<C(y_0)=0,\qquad \forall\ k>0,$$
so that  
$$g'y_{k}<-y_{k}'Hy_{k}<0,\qquad \forall\ k>0.\xdef\descent{\lab}\eqnum\show{oneo}$$
We may either let the process terminate naturally (which will happen after a number of iterations no larger than the number of paths $P$), or more practically, terminate once a certain termination criterion is satisfied, in which case $y_k$ is a descent direction by Eq.\ \descent.

An important variant, which is usually far superior in practice, is to use a preconditioning matrix $S$, which is a diagonal approximation to the Hessian matrix, i.e, $S$ is diagonal with $H_{pp}$ along the diagonal. This method starts with $y_0=0$, $r_0=g$, and $p_0=-Sg$, and has the form
$$\eqalign{y_{k+1}&=y_k+\a_kp_k,\qquad\hbox{where }\a_k={r_k'Sr_k\over p_k'Hp_k},\cr
r_{k+1}&=g+Hy_{k+1},\cr
p_{k+1}&=-Sr_{k+1}+\b_kp_k,\qquad\hbox{where }\b_k={r_{k+1}'Sr_{k+1}\over r_k'Sr_k}.\cr}$$
The non-preconditioned method bears the same relation to the preconditioned version that the gradient method bears to the diagonally scaled gradient method, with second-order diagonal scaling. The use of this type of preconditioning not only improves the rate of convergence (typically), but also facilitates the choice of stepsize (a stepsize of 1  typically works, regardless of the number of conjugate gradient iterations used). These advantages of preconditioning have been confirmed by extensive computational experience, although theoretically speaking there are rare exceptions where diagonal second-order preconditioning does not improve performance. 

A somewhat different type of diagonal preconditioning can be used with advantage in the case where the number of arcs, call it $A$, is substantially smaller than the number of paths $P$. Then by using $S=\gr^2R(x)$ as a preconditioning matrix, it can be shown that the number of conjugate gradient iterations to find the Newton direction is at most $A$, as noted in [Ber74] (see also [Ber99], p.\ 148).

In the preceding methods, each conjugate gradient iteration involves a small number of inner products and matrix-vector multiplications, each of which requires no more than $O(T)$ computation. Thus the total required number of computational operations to compute approximately the Newton direction with $m$ conjugate gradient iterations is $O(mT)$, compared to $O(T)$ for a single gradient calculation and diagonally scaled gradient iteration.

The choice of number of conjugate gradient steps to obtain an approximate Newton direction is an interesting practical implementation issue, and may ultimately be settled by experimentation. The paper [DES82] shows how to control this number so that a linear, superlinear, or quadratic convergence rate for the overall method is achieved. Generally, to attain a superlinear rate, the conjugate gradient process for approximating the Newton step must become asymptotically exact, with the number of conjugate gradient steps per iteration approaching $P$. This appears highly inefficient for practical large network flow problems, where $P$ is large and just a few conjugate gradient steps  per iteration are sufficient to attain a good convergence rate, particularly when second order diagonal preconditioning is used.

\vskip-1.5pc

\section{Distributed Computation of the Newton Step}
\vskip-1.0pc
\pn The algorithm of the preceding section can also be executed in a distributed fashion, assuming that there is a processor assigned to each path and a processor assigned to each arc. We may also assign multiple paths and arcs to a single processor, with an attendant increase of the computation and communication load of the processor, but for simplicity we do not consider this possibility. 

To compute the Newton step, for each $p$, the processor assigned to path $p$ may update, maintain, and communicate to other processors its own path variable $x_p$, and compute the corresponding component $y_p$ of the Newton direction. For each $a$, the processor assigned to arc $a$ may compute and communicate to the relevant path processors, its accumulated flow variable $f_a$. Similar computations are used to execute the intermediate conjugate gradient iterations. 
The path processors may also collaborate to compute a stepsize guaranteeing descent at every iteration, although in many network applications a constant stepsize can be used reliably, as experience has shown [BGG84], [GaB84] (equal or nearly equal to one in our case). 

Information exchange between the processors may be conveniently and accurately performed by using a spanning tree with a designated node serving as a synchronizer and leader for the distributed computation. Schemes of this type have been used extensively in data networks and distributed computing systems. Recently, iterative consensus schemes [TBA86], [BeT89] have been discussed as possible methods for information exchange between processors, but these schemes are slow relative to Newton's method, so they are not suitable for our context. On the other hand, it should be noted that in certain application contexts, the need for synchronization to perform Newton and conjugate gradient iterations may be a major drawback over asynchronous diagonally scaled gradient and gradient projection methods, which can be implemented asynchronously with satisfactory convergence properties (see [Ber83], [BeT89] for totally asynchronous gradient-like methods, and [TsB86], [BeT89] for partially asynchronous and stochastic gradient-like methods).

\section{Algorithms, Convergence, and Rate of Convergence}

\vskip-1.0pc

\pn There are several network flow problem formulations and corresponding algorithms, where the Newton direction computations of the preceding two sections may be potentially used.
An important issue is the treatment of constraints in this context. Here are some potential approaches, leading to viable Newton-type algorithms:

\nitem{(a)} Constraints may be eliminated via a penalty or augmented Lagrangian approach to yield an unconstrained problem of minimizing a function $F(x)$ of the form \costfn. The convergence analysis of such approaches requires (1) a convergence guarantee for the unconstrained optimization of the penalized or augmented Lagrangian objective, and (2) a convergence guarantee for the overall penalty or augmented Lagrangian objective. These guarantees may be obtained in straightforward fashion as applications of standard results for gradient-related and Newton-like methods for unconstrained minimization, as well as standard analyses of penalty and augmented Lagrangian methods (see e.g., [Ber82a], [Ber99]).

\nitem{(b)} Nonnegativity constraints can be treated with two-metric Newton-like methods that require computation of Newton directions such as the ones of the preceding two sections. These are methods involving partial diagonalization of the Hessian matrix, for which a detailed convergence and rate of convergence analysis is given in [Ber82b], [BeG83], [GaB84].

\nitem{(c)} Simplex constraints, relating to supply/demand specifications, can be used to eliminate some of the variables and essentially reduce the constraint set to the nonnegativity case of (b) above (see again [Ber82b], [BeG83], [GaB84] for convergence and rate of convergence analysis).

\smskip

Regardless of how constraints are treated, for a convergent algorithm one must address the issue of stepsize selection. Experience has shown that in network algorithms one may often use reliably a constant stepsize, particularly when the directions used embody second order information, which makes a stepsize close to one a typically good choice. An alternative to a constant stepsize is a line search rule based on line minimization or successive stepsize reduction. References [Ber82b], [BeG83], [GaB84] provide examples of successive stepsize reduction rules in conjunction with constraints. The use of such rules improves the reliability of algorithms but introduces additional complexity, particularly in a distributed context. In the latter case, it is possible to implement successive stepsize reduction rules at a special processor that may exchange information with other processors in a distributed way.

We do not discuss convergence issues further. The aim of this paper is not to provide specific algorithms and associated convergence analysis (which is routine for the most part, as well as problem dependent), but rather to make the point that the network structure can be used to implement the computation of Newton directions in a convenient centralized or distributed manner.

\section{Extensions}

\vskip-0.5pc

\pn There are several problem generalizations, involving in some cases an extended network or even a non-network structure, which admit a treatment  similar to the one of the preceding sections:

\nitem{(a)} Cost functions $F$ that are defined over just a subset of the path flow space arise when constraints are eliminated by means of an interior point method approach. The use of Newton directions within this context is well-documented.

\nitem{(b)} More general linear dependence of arc flows on path flows 
 can be treated by generalization of the terms $D_a(f_a)$ in the cost function. In particular, we may redefine $f_a$ to be a general linear function $c_a'x$ where $c_a$ is some vector. For example the scalar components of $c_a$ may represent arc gains. In this case, the entire approach of the preceding two sections generalizes straightforwardly. What is essential is that the Hessian of $F$ should have the generic form 
 $$\gr^2 R+E'\gr^2D E$$
  of Eq.\ \hescomp, with $\gr^2 R$ and $\gr^2D$ being diagonal matrices, and $E$ being a matrix that encodes a linear dependence between $x$ and the arguments $f_a$ of the cost terms $D_a(f_a)$.

\nitem{(c)} Linear side constraints may be treated by using a penalty or augmented Lagrangian approach, thereby reducing to case (b) above.

\nitem{(d)} The basic Hessian structure that is important for the convenient computation of gradients and Hessian matrix-vector products is 
$$E'DE,$$
where $D$ and $E$ are matrices with $D$ diagonal. Therefore our methodology will also work for structures of the form
$$R+E_1'D_1E_1+\cdots+E_m'D_mE_m,$$
where $R$, $D_1,\ldots,D_m$, and $E_1,\ldots,E_m$ are matrices with $R$ and $D_1,\ldots,D_m$ diagonal. There are also potential extensions in cases where $R$ and $D_1,\ldots,D_m$ are symmetric and nearly diagonal (e.g., tridiagonal).

\smskip

An important question is how to deal with singularity of the Hessian matrix and the attendant lack of strong convexity. This arises for example in the important case where $\gr^2 R=0$, and there are constraints on $x$ (nonnegativity and/or supply/demand constraints), which guarantee existence of a solution. While in this case subsequence convergence  of our methods to minima is easy to show under standard assumptions, the convergence (to a single point) and the establishment of a linear or superlinear convergence rate result are open questions. By contrast, the issue of convergence  to a single point and linear convergence rate in the presence of Hessian singularity has been satisfactorily addressed for gradient projection methods in the context of variational inequalities, including multicommodity flow problems [BeG82].

\medskip\smskip

\pn {\bf \hskip0pc REFERENCES}

\medskip

\def\ref{\vskip1pt\pn}

\ref[AtL00] Athuraliya, S., and Low, S., 2000.\ ``Optimization Flow Control with Newton-Like Algorithm," J.\ of Telecommunication Systems, Vol.\ 15, pp.\ 345-358.

\ref[BGG84] Bertsekas, D.\ P., Gafni, E.\ M., and Gallager, R.\ G., 1984.\ ``Second
Derivative Algorithms for Minimum Delay Distributed Routing in
Networks," IEEE Trans.\ on Communications, Vol.\ 32, pp.\ 911-919.

\ref[BeG82]  Bertsekas, D.\ P., and Gafni, E.\ M., 1982.\ ``Projection Methods for
Variational Inequalities with Application to the Traffic Assignment
Problem," Math.\ Progr.\ Studies, Vol.\ 17, North-Holland, Amsterdam, pp.\
139-159.

\ref[BeG83] Bertsekas, D.\ P., and Gafni, E.\ M., 1983.\ ``Projected Newton Methods
and Optimization of Multicommodity Flows," IEEE Trans.\ on Auto.
Control, Vol.\ 28, pp.\ 1090-1096.

\ref[BeG92] Bertsekas, D.\ P., and Gallager, R.\ G., 1992.\ Data
Networks, (2nd Ed.), Prentice-Hall, Englewood Cliffs, N.\ J.

\ref [BeT89] Bertsekas, D.\ P., and Tsitsiklis, J.\ N., 1989.\ Parallel and
Distributed Computation: Numerical Methods, Prentice-Hall, Englewood Cliffs,
N.\ J.\ (republished in 1997 by Athena Scientific, Belmont, MA).

\ref[BeY10] Bertsekas, D.\ P., and Yu, H., 2010.\ ``A Unifying Polyhedral Approximation Framework for Convex Optimization," Lab.\ for Information and Decision Systems Report LIDS-P-2820, MIT; to appear in SIAM J. on Optimization.

\ref [Ber74] Bertsekas, D.\ P., 1974.\  ``Partial Conjugate Gradient Methods
for a Class of Optimal Control Problems," IEEE Trans.\ Automat.\ Control, Vol.\
19, pp.\ 209-217.

\ref [Ber80] Bertsekas, D.\ P., 1980.\  ``A Class of Optimal Routing Algorithms for
Communication Networks," Proc.\ of the Fifth International Conference on
Computer Communication, Atlanta, GA, Oct. 1980, pp.\ 71-76.

\ref [Ber82a] Bertsekas, D.\ P., 1982.\ Constrained Optimization and Lagrange
Multiplier Methods, Academic Press, N.\ Y.\ 
(republished in 1996 by Athena Scientific, Belmont, MA).

\ref [Ber82b] Bertsekas, D.\ P., 1982.\  ``Projected Newton Methods for
Optimization Problems with Simple Constraints," SIAM J.\ on Control and Optimization, Vol.\
20, pp.\ 221-246.

\ref[Ber83] Bertsekas, D.\ P.,  1983.\ ``Distributed Asynchronous Computation of Fixed Points," Mathematical Programming, Vol.\ 27, pp.\ 107-120.

\ref[Ber98] Bertsekas, D.\ P., 1998.\ Network Optimization: Continuous and Discrete Models, Athena Scientific, Belmont, MA. 

\ref[Ber99] Bertsekas, D.\ P., 1999.\ Nonlinear Programming, Athena Scientific, Belmont, MA. 

\ref [CaG74] Cantor, D.\ G., and Gerla, M., 1974.\ ``Optimal Routing in a Packet
Switched Computer Network," IEEE Trans.\ on Computers, Vol.\ 23, pp.\
1062-1069.

\ref [DES82] Dembo, R.\ S., Eisenstadt, S.\ C., and Steihaug, T., 1982.\ 
``Inexact Newton Methods," SIAM J.\ Numer.\ Anal., Vol.\ 19, pp.\ 400-408.

\ref[FGK73] Fratta, L., Gerla, M., and Kleinrock, L., 1973.\ ``The Flow-Deviation Method: An
Approach to Store-and-Forward Computer Communication Network Design,''
Networks, Vol.\ 3, pp.\ 97-133.

\ref[FlH95] Florian, M.\ S., and Hearn, D., 1995.\ ``Network Equilibrium Models and Algorithms,'' 
Handbooks in OR and MS,  Ball, M.\ O., Magnanti, T.\ L., Monma, C.\ L., and Nemhauser, G.\ L.\
(eds.),  Vol.\ 8, North-Holland, Amsterdam, pp.\ 485-550. 

\ref [GaB84] Gafni, E.\ M., and Bertsekas, D.\ P., 1984.\  ``Two-Metric Projection Methods for
Constrained Optimization," SIAM J.\ on Control and Optimization, Vol.\ 22, pp.\ 936-964.

\ref[GaG80] Gallager, R.\ G., and Golestaani, S.\ J., 1980.\ ``Flow Control and Routing Algorithms for Data Networks," Proc.\ 5th Intern.\ Conf.\ Comput.\ Comm., pp.\ 779-784.

\ref[Gal77] Gallager, R.\ G., 1977.\ ``A Minimum Delay Routing Algorithm Using
Distributed Computation," IEEE Trans.\ on Communications, Vol.\ 23, pp.\
73-85.

\ref[HLV87] Hearn, D.\ W., Lawphongpanich, S., and Ventura, J.\ A., 1987.\ ``Restricted
Simplicial Decomposition: Computation and Extensions," Math.\ Programming Studies, Vol.\
31, pp.\ 119-136.

\ref[Hoh77] Hohenbalken, B.\ von, 1977.\ ``Simplicial Decomposition in Nonlinear
Programming," Math.\ Programming, Vol.\ 13, pp.\ 49-68.

\ref[Hol74] Holloway, C.\ A., 1974.\ ``An Extension of the Frank and Wolfe Method of
Feasible Directions," Math.\ Programming, Vol.\ 6, pp.\ 14-27.

\ref[JOZ09] Jadbabaie, A., Ozdaglar, A., and Zargham, M., 2009.\ ``A Distributed Newton Method for Network Optimization," Proc.\ of 2009 CDC.

\ref [Lue84] Luenberger, D.\ G., 1984.\  Introduction to Linear and Nonlinear
Programming, (2nd Ed.), Addison-Wesley, Reading, MA.

\ref [Pat98] Patriksson, M., 1998.\ Nonlinear Programming and Variational Inequalities: A Unified
Approach, Kluwer, Dordtrecht, The Netherlands.

\ref [Roc84] Rockafellar, R.\ T., 1984.\ Network Flows and Monotropic
Programming, Wiley, N.\ Y.

\ref[TBA86] Tsitsiklis, J.\ N., Bertsekas, D.\ P., and Athans, M., 1986.\
``Distributed Asynchronous Deterministic and Stochastic Gradient Optimization Algorithms," IEEE Trans.\
on Aut.\ Control, Vol.\ AC-31, pp.\ 803-812.

\ref [TsB86] Tsitsiklis, J.\ N., and Bertsekas, D.\ P., 1986.\ ``Distributed Asynchronous Optimal
Routing in Data Networks,'' IEEE Trans.\ on Automatic Control, Vol.\ 31, pp.\ 325-331.

\ref[WOJ10] Wei, E., Ozdaglar, A., and Jadbabaie, A., 2010.\ ``A Distributed Newton Method for Network Utility Maximization,"Lab.\ for Information and Decision Systems Report LIDS-2832, M.I.T.; also in Proc.\ of 2010 CDC.

\end

%% file: TEXSHOP_macros_new.tex

\def\ignore#1{}
 

\newcount\sectnum
\newcount\subsectnum
\newcount\eqnumber

\global\eqnumber=1\sectnum=0


\def\lab{(\the\sectnum.\the\eqnumber)}



\def\show#1{#1}



\def\smskip{\vskip 5 pt}
\def\medskip{\vskip 10 pt}
\def\bigskip{\vskip 15 pt}
\def\pn{\par\noindent}

\def\frac#1#2{{#1\over #2}}

\def\half{{\scriptstyle {1\over 2}}}

\def\a{\alpha}

\def\b{\beta}

\def\re{\Re}

\def\gr{\nabla}


\def\old#1{}
\def\leaderfill{\leaders\hbox to 1em{\hss.\hss}\hfill}


\parindent=2pc
\baselineskip=15pt
\vsize=8.7 true in
\voffset=0.125 true in
\parskip=3pt


\def\minprob#1#2#3{$$\eqalign{&\hbox{minimize\ \ }#1\cr &\hbox{subject to\ \
}#2\cr}\ifnum 0=#3{}\else\eqno(#3)\fi$$}        
     
\def\maxprob#1#2#3{$$\eqalign{&\hbox{maximize\ \ }#1\cr &\hbox{subject to\ \
}#2\cr}\ifnum 0=#3{}\else\eqno(#3)\fi$$}        
     
\def\aligntwo#1#2#3#4#5{$$\eqalign{#1&#2\cr #3&#4\cr}
\ifnum 0=#5{}\else\eqno(#5)\fi$$}
\def\alignthree#1#2#3#4#5#6#7{$$\eqalign{#1&#2\cr #3&#4\cr #5&#6\cr}
\ifnum 0=#7{}\else\eqno(#7)\fi$$}


\def\eqnum{\eqno{\hbox{(\the\sectnum.\the\eqnumber)}\global\advance\eqnumber
by1}}

\def\eqnu{\eqno{\hbox{(\the\sectnum.\the\eqnumber)}\global\advance\eqnumber
by1}}

\newcount\examplnumber
\def\examplnum{\global\advance\examplnumber by1}

\newcount\figrnumber
\def\figrnum{\global\advance\figrnumber by1}

\newcount\propnumber
\def\propnum{\global\advance\propnumber by1}

\newcount\defnumber
\def\defnum{\global\advance\defnumber by1}

\newcount\lemmanumber
\def\lemmanum{\global\advance\lemmanumber by1}

\newcount\assumptionnumber
\def\assumptionnum{\global\advance\assumptionnumber by1}

\newcount\conditionnumber
\def\conditionnum{\global\advance\conditionnumber by1}

\def\exampl{\the\sectnum.\the\examplnumber}
\def\figr{\the\sectnum.\the\figrnumber}
\def\propn{\the\sectnum.\the\propnumber}
\def\defn{\the\sectnum.\the\defnumber}
\def\lemman{\the\sectnum.\the\lemmanumber}
\def\assumptionn{\the\sectnum.\the\assumptionnumber}
\def\condn{\the\sectnum.\the\conditionnumber}

\def\section#1{\goodbreak\vskip 3pc plus 6pt minus 3pt\leftskip=-2pc
   \global\advance\sectnum by 1\eqnumber=1
\global\examplnumber=1\figrnumber=1\propnumber=1\defnumber=1\lemmanumber=1\assumptionnumber=1 \conditionnumber =1%
   \line{\hfuzz=1pc{\hbox to 3pc{\bf 
   \vtop{\hfuzz=1pc\hsize=38pc\hyphenpenalty=10000\noindent\uppercase{\the\sectnum.\quad #1}}\hss}}
			\hfill}
			\leftskip=0pc\nobreak\tenf
			\vskip 1pc plus 4pt minus 2pt\noindent\ignorespaces}



\def\sect#1{\noindent\leftskip=-2pc\tenf
   \goodbreak\vskip 1pc plus 4pt minus 2pt
                \global\advance\subsectnum by 1\eqnumber=1
   \line{\hfuzz=1pc{\hbox to 3pc{\bf 
   \vtop{\hfuzz=1pc\hsize=38pc\hyphenpenalty=10000\noindent\uppercase{{\bf #1}}}\hss}}
                        \hfill}
   \leftskip=0pc\nobreak\tenf
                        \vskip 1pc plus 4pt minus 2pt\nobreak\noindent\ignorespaces}

\def\subsection#1{\noindent\leftskip=0pc\tenf
   \goodbreak\vskip 1pc plus 4pt minus 2pt
   \line{\hfuzz=1pc{\hbox to 3pc{\bf 
   \vtop{\hfuzz=1pc\hsize=38pc\hyphenpenalty=10000\noindent{\bf #1}}\hss}}
                        \hfill}
   \leftskip=0pc\nobreak\tenf
                        \vskip 1pc plus 4pt minus 2pt\nobreak\noindent\ignorespaces}
\def\subsubsection#1{\goodbreak\vskip 1pc plus 4pt minus 2pt
   \hfuzz=3pc\leftskip=0pc\noindent\tenit #1 \nobreak\tenf\vskip 6pt plus 1pt
                                minus 1pt\nobreak\ignorespaces\leftskip=0pc}
%

\def\beginexample#1{\noindent\goodbreak\vskip 6pt plus 1pt minus 1pt
\noindent
  \hbox {\bf Example #1\hss}
  \nobreak\vskip 4pt plus 1pt minus 1pt \nobreak\noindent\ninef
  \global\advance
                \leftskip by\parindent\pn}
\def\endexample{\vskip 12pt\tenf\par
  \global\advance\leftskip by -\parindent
  }

\def\beginexercise#1{\noindent\goodbreak\vskip 6pt plus 1pt minus 1pt \noindent\global\normalbaselineskip=12pt
  \hbox {\bf Exercise #1\hss}
  \nobreak\vskip 4pt plus 1pt minus 1pt 
  \nobreak\noindent\ninef\global\advance\leftskip
                        by\parindent\pn}
\def\endexercise{\vskip 12pt\tenf\par
  \global\advance\leftskip by -\parindent
  }

\def\beginsection#1{\noindent\goodbreak\vskip 6pt plus 1pt minus 1pt \noindent\global\normalbaselineskip=12pt
  \hbox {\it #1\hss}
  \vskip 0.1pt plus 1pt minus 1pt \nobreak\noindent\ninef\global\advance
                \leftskip by\parindent\noindent\pn}
\def\endsection{\vskip 12pt\tenf\par
  \global\advance\leftskip by -\parindent
}

%


\def\ref{\smskip\pn}

\def\chapter#1#2{{\bf \centerline{\helbigbig
{#1}}}\bigskip\bigskip{\bf \centerline{\helbigbig
{#2}}}\bigskip\bigskip} 



\def\longpapertitle#1#2#3{{\bf \centerline{\helbigb
{#1}}}\bigskip{\bf \centerline{\helbigb
{#2}}}\bigskip\bigskip{\centerline{
by}}\bigskip{\bf \centerline{
{#3}}}\bigskip\bigskip} 


\def\nitem#1{\smskip\item{#1}}

\newcount\alphanum
\newcount\romnum

\def\alphaenumerate{\ifcase\alphanum \or (a)\or (b)\or (c)\or (d)\or (e)\or
(f)\or (g)\or (h)\or (i)\or (j)\or (k)\fi}
\def\romenumerate{\ifcase\romnum \or (i)\or (ii)\or (iii)\or (iv)\or (v)\or
(vi)\or (vii)\or (viii)\or (ix)\or (x)\or (xi)\fi}

\def\alist{\begingroup\vskip10pt\alphanum=1
\parskip=2pt\parindent=0pt \leftskip=3pc
\everypar{\llap{{\rm\alphaenumerate\hskip1em}}\advance\alphanum by1}}

\def\nolist{\begingroup\vskip10pt\alphanum=0
\parskip=2pt\parindent=0pt \leftskip=3pc
\everypar{\llap{\global\advance\alphanum by1(\the\alphanum)\hskip1em}}}

\def\romlist{\begingroup\vskip10pt\romnum=1
\parskip=2pt\parindent=0pt \leftskip=5pc
\everypar{\llap{{\rm\romenumerate\hskip1em}}\advance\romnum by1}}



\long\def\fig#1#2#3{\vbox{\vskip1pc\vskip#1
\prevdepth=12pt \baselineskip=12pt
\vskip1pc
\hbox to\hsize{\hfill\vtop{\hsize=25pc\noindent{\eightbf Figure #2\ }
{\eightpoint#3}}\hfill}}}

\long\def\widefig#1#2#3{\vbox{\vskip1pc\vskip#1
\prevdepth=12pt \baselineskip=12pt
\vskip1pc
\hbox to\hsize{\hfill\vtop{\hsize=28pc\noindent{\eightbf Figure #2\ }
{\eightpoint#3}}\hfill}}}

\long\def\table#1#2{\vbox{\vskip0.5pc
\prevdepth=12pt \baselineskip=12pt
\hbox to\hsize{\hfill\vtop{\hsize=25pc\noindent{\eightbf Table #1\ }
{\eightpoint#2}}\hfill}}}

 
\def\rightheadline#1{\headline{\tenrm\hfil #1}}


\long\def\leftfig#1#2{\vbox{\smskip\hsize=220pt
\vtop{{\noindent {\bf #1}}}
\smskip
\noindent
\vbox{{\noindent #2}}
}}

\long\def\rightfig#1#2#3{\vbox{\smskip\vskip#1
\prevdepth=12pt \baselineskip=12pt
\hsize=210pt
\smskip
\vbox{\noindent{\eightbold #2}
\hskip1em{\eightpoint#3}}
}}

\long\def\concept#1#2#3#4#5{\bigskip\hrule
\vbox{\hbox{\leftfig{#1}{#2} \hskip3em
\rightfig{#3}{#4}{#5}} \smskip}
\hrule\bigskip}


\long\def\bconcept#1#2#3#4#5#6#7{
\vbox{
\hbox to \hsize{\vtop{\par #1}}
\concept{#2}{#3}{#4}{#5}{#6}
\hbox to \hsize{\vtop{\par #7}}
\smskip}
}




\def\boxit#1{\vbox{\hrule\hbox{\vrule\kern3pt
                                \vbox{\kern3pt#1\kern3pt}\kern3pt\vrule}\hrule}}
\def\centerboxit#1{$$\vbox{\hrule\hbox{\vrule\kern3pt
                                \vbox{\kern3pt#1\kern3pt}\kern3pt\vrule}\hrule}$$}


%
%
%

\def\picture #1 by #2 (#3){
  \vbox to #2{
    \hrule width #1 height 0pt depth 0pt
    \vfill
    \special{picture #3} 
    }
  }

\def\scaledpicture #1 by #2 (#3 scaled #4){{
  \dimen0=#1 \dimen1=#2
  \divide\dimen0 by 1000 \multiply\dimen0 by #4
  \divide\dimen1 by 1000 \multiply\dimen1 by #4
  \picture \dimen0 by \dimen1 (#3 scaled #4)}
  }

%
%

\long\def\captfig#1#2#3#4#5{\vbox{\vskip1pc
\hbox to\hsize{\hfill{\picture #1 by #2 (#3)}\hfill}
\prevdepth=9pt \baselineskip=9pt
\vskip1pc
\hbox to\hsize{\hfill\vtop{\hsize=24pc\noindent{\eightbold Figure #4}
\hskip1em{\eightpoint#5}}\hfill}}}

%
%
%

\def\illustration #1 by #2 (#3){
  \vskip#2\hskip#1\special{illustration #3} 
    }

\def\scaledillustration #1 by #2 (#3 scaled #4){{
  \dimen0=#1 \dimen1=#2
  \divide\dimen0 by 1000 \multiply\dimen0 by #4
  \divide\dimen1 by 1000 \multiply\dimen1 by #4
  \illustration \dimen0 by \dimen1 (#3 scaled #4)}
  }


\newbox\graybox
\newdimen\xgrayspace
\newdimen\ygrayspace
%
%
%
%
%
%
%
%
%

\def\Textshade#1#2#3#4#5#6{%
    \xgrayspace=#4pt%
    \ygrayspace=#4pt%
    \def\grayshade{#3}%
    \def\linewidth{#5}%
    \def\theradius{#6}%
    \setbox\graybox=\hbox{\surroundboxa{#2}}%
    \hbox{%
    \hbox to 0pt{%
    \PScommands
}%
    \box\graybox}}%
%
%
\long%

\long%
\def\Parashade#1#2#3#4#5#6#7{%
    \xgrayspace=#4pt%
    \ygrayspace=#4pt%
    \def\grayshade{#3}%
    \def\linewidth{#5}%
    \def\theradius{#6}%
    \def\thevskip{#7pt}%
    \setbox\graybox=\hbox{\surroundboxb{#2}}%
    \vskip\thevskip%
    \hbox{%
    \hbox to 0pt{%
    \PScommands
}%
     \box\graybox}%
     \vskip\thevskip%
}%
%
%
%
\long\def\surroundboxa#1{\leavevmode\hbox{\vtop{%
\vbox{\kern\ygrayspace%
\hbox{\kern\xgrayspace#1%
      \kern\xgrayspace}}\kern\ygrayspace}}}
%
%
\long\def\surroundboxb#1{\leavevmode\hbox{\vtop{%
\vbox{\kern\ygrayspace%
\hbox{\kern\xgrayspace\vbox{\advance\hsize-2\xgrayspace#1}%
      \kern\xgrayspace}}\kern\ygrayspace}}}
%
%
%
\long\def\PScommands{%
\special{rawpostscript
/sharpbox{%
           newpath
           xmin ymin moveto
           xmin ymax lineto
           xmax ymax lineto
           xmax ymin lineto
           xmin ymin lineto
           closepath 
          }bind def
}%
\special{rawpostscript
/sharpboxnb{%
           newpath
           xmin ymin moveto
           xmin ymax lineto
           xmax ymax lineto
           xmax ymin lineto
          }bind def
}%
\special{rawpostscript
/sharpboxnt{%
           newpath
           xmin ymax moveto
           xmin ymin lineto
           xmax ymin lineto
           xmax ymax lineto
          }bind def
}%
\special{rawpostscript
/roundbox{%
           newpath
           xmin radius add ymin moveto
           xmax ymin xmax ymax radius arcto
           xmax ymax xmin ymax radius arcto
           xmin ymax xmin ymin radius arcto
           xmin ymin xmax ymin radius arcto 16 {pop} repeat
           closepath
          }bind def
}%
\special{rawpostscript
/sharpcorners{%
               sharpbox gsave grayshade setgray fill grestore 
               linewidth setlinewidth stroke
              }bind def
}%
\special{rawpostscript
/sharpcornersnt{%
               sharpboxnt gsave grayshade setgray fill grestore 
               linewidth setlinewidth stroke
              }bind def
}%
\special{rawpostscript
/sharpcornersnb{%
               sharpboxnb gsave grayshade setgray fill grestore 
               linewidth setlinewidth stroke
              }bind def
}%
\special{rawpostscript
/roundcorners{%
               roundbox gsave grayshade setgray fill grestore 
               linewidth setlinewidth stroke
              }bind def
}%
\special{rawpostscript
/plainbox{%
           sharpbox grayshade setgray fill 
          }bind def
}%
%
\special{rawpostscript
/roundnoframe{%
               roundbox grayshade setgray fill 
              }bind def
}%
\special{rawpostscript
/sharpnoframe{%
               sharpbox grayshade setgray fill 
              }bind def
}%
}%
%
%


\def\boxit#1{\vbox{\hrule\hbox{\vrule\kern3pt
                                \vbox{\kern3pt#1\kern3pt}\kern3pt\vrule}\hrule}}

\def\boxitnb#1{\vbox{\hrule\hbox{\vrule\kern3pt
                                \vbox{\kern3pt#1\kern3pt}\kern3pt\vrule}}}

\def\boxitnt#1{\vbox{\hbox{\vrule\kern3pt
                                \vbox{\kern3pt#1\kern3pt}\kern3pt\vrule}\hrule}}





%
%
%
%
%
%
%
%
\font\helbigbig=cmr10 scaled 2500%
\font\helbigb=cmbx10 scaled 1500%
\font\eightbold=cmbx8%

\def\tenf{\hel}%
\def\tenit{\heli}%
\def\ninef{\ninehel}%
\def\nineit{\nineheli}%
%
%


\font\tenrm=cmr10%
\font\teni=cmmi10%
\font\tensy=cmsy10%
\font\tenbf=cmbx10%
\font\tentt=cmtt10%
\font\tenit=cmti10%
\font\tensl=cmsl10%

\def\tenpoint{\def\rm{\fam0\tenrm}%
\textfont0=\tenrm%
\textfont1=\teni%
\textfont2=\tensy%
\textfont\itfam=\tenit%
\textfont\slfam=\tensl%
\textfont\ttfam=\tentt%
\textfont\bffam=\tenbf%
\scriptfont0=\sevenrm%
\scriptfont1=\seveni%
\scriptfont2=\sevensy%
\scriptscriptfont0=\sixrm%
\scriptscriptfont1=\sixi%
\scriptscriptfont2=\sixsy%
\def\it{\fam\itfam\tenit}%
\def\tt{\fam\ttfam\tentt}%
\def\sl{\fam\slfam\tensl}%
\scriptfont\bffam=\sevenbf%
\scriptscriptfont\bffam=\sixbf%
\def\bf{\fam\bffam\tenbf}%
\normalbaselineskip=18pt%
\normalbaselines\rm}%

\font\ninerm=cmr9%
\font\ninebf=cmbx9%
\font\nineit=cmti9%
\font\ninesy=cmsy9%
\font\ninei=cmmi9%
\font\ninett=cmtt9%
\font\ninesl=cmsl9%

\def\ninepoint{\def\rm{\fam0\ninerm}%
\textfont0=\ninerm%
\textfont1=\ninei%
\textfont2=\ninesy%
\textfont\itfam=\nineit%
\textfont\slfam=\ninesl%
\textfont\ttfam=\ninett%
\textfont\bffam=\ninebf%
\scriptfont0=\sixrm%
\scriptfont1=\sixi%
\scriptfont2=\sixsy%
\def\it{\fam\itfam\nineit}%
\def\tt{\fam\ttfam\ninett}%
\def\sl{\fam\slfam\ninesl}%
\scriptfont\bffam=\sixbf%
\scriptscriptfont\bffam=\fivebf%
\def\bf{\fam\bffam\ninebf}%
\normalbaselineskip=16pt%
\normalbaselines\rm}%

\font\eightrm=cmr8%
\font\eighti=cmmi8%
\font\eightsy=cmsy8%
\font\eightbf=cmbx8%
\font\eighttt=cmtt8%
\font\eightit=cmti8%
\font\eightsl=cmsl8%

\def\eightpoint{\def\rm{\fam0\eightrm}%
\textfont0=\eightrm%
\textfont1=\eighti%
\textfont2=\eightsy%
\textfont\itfam=\eightit%
\textfont\slfam=\eightsl%
\textfont\ttfam=\eighttt%
\textfont\bffam=\eightbf%
\scriptfont0=\sixrm%
\scriptfont1=\sixi%
\scriptfont2=\sixsy%
\scriptscriptfont0=\fiverm%
\scriptscriptfont1=\fivei%
\scriptscriptfont2=\fivesy%
\def\it{\fam\itfam\eightit}%
\def\tt{\fam\ttfam\eighttt}%
\def\sl{\fam\slfam\eightsl}%
\scriptscriptfont\bffam=\fivebf%
\def\bf{\fam\bffam\eightbf}%
\normalbaselineskip=14pt%
\normalbaselines\rm}%

\font\sevenrm=cmr7%
\font\seveni=cmmi7%
\font\sevensy=cmsy7%
\font\sevenbf=cmbx7%

\font\sixrm=cmr6%
\font\sixi=cmmi6%
\font\sixsy=cmsy6%
\font\sixbf=cmbx6%

\fontdimen13\tensy=2.6pt%
\fontdimen14\tensy=2.6pt%
\fontdimen15\tensy=2.6pt%
\fontdimen16\tensy=1.2pt%
\fontdimen17\tensy=1.2pt%
\fontdimen18\tensy=1.2pt%

\def\tenf{\tenpoint}%
\def\ninef{\ninepoint}%
%
